\documentclass{amsart}
\usepackage{amssymb, amsmath, latexsym}

\renewcommand{\baselinestretch}{\baselinestretch}
\renewcommand{\baselinestretch}{1.1}
\author{Constantin-Nicolae Beli}
\title[]{Two results on cohomology of groups adapted to cochains} 
\date{}

\def\m{\lim}

 \def\({\overline}

\def\){\underline} \def\<{\cdot} 
\def\>{~~~~~~~} \def\#{{\bf
Definition}} \def\*{\section} \def\be{\begin{equation}}
\def\ee{\end{equation}}

\def\bmat{\left(\begin{array}} \def\emat{\end{array}\right)}

 \def\m2{~(\mo 2)} \def\no{\noindent}
 \def\btm{\begin{thm}}
\def\etm{\end{tm}}
 \def\blem{\begin{lem}}
\def\elem{\end{lem}}

\newtheorem{theorem}{Theorem}[section]
\newtheorem{proposition}[theorem]{Proposition}
\newtheorem{lemma}[theorem]{Lemma}
\newtheorem{definition}{Definition}
\newtheorem{corollary}[theorem]{Corollary}

\newtheorem{bof}[theorem]{}
\newtheorem{teorema}{Theorem}

\def\qed{\mbox{$\Box$}\vspace{\baselineskip}}
\def\pf{$Proof.$ } 
\def\bco{\begin{corollary}} \def\eco{\end{corollary}} 
\def\bdf{\begin{definition}} \def\edf{\end{definition}} 
\def\btm{\begin{theorem}} \def\etm{\end{theorem}} 
\def\bpr{\begin{proposition}} \def\epr{\end{proposition}}  
\def\blm{\begin{lemma}} \def\elm{\end{lemma}} 
\def\bff{\begin{bof}\rm} \def\eff{\end{bof}}
\def\btr{\begin{teorema}} \def\etr{\end{teorema}}

\def\de{\newcommand} \de\tm[1]{{\no\bf Theorem~#1}} 
\def\mb{\mathbb} 
   \def\ZZ{{\mb Z}}
 
\de\lm[1]{{\no\bf Lemma~#1}}
\de\df[1]{{\no\bf Definition~#1}} \de\co[1]{{\no\bf Corollary~#1}}
\de\lr[1]{\longrightarrow^{\!\!\!\!\!\!\!\! #1}}
\de\lf[1]{\longleftarrow^{\!\!\!\!\!\!\!\! #1}}
\de\si[1]{\sim^{\!\!\!\!\! #1}} \de\apr[1]{\approx^{\!\!\!\!\! #1}}
\de\leg[2]{\left(\frac {#1}{#2}\right)}
\DeclareMathOperator\Br{Br}

\de\Brr[1]{{}_{#1}\Br}

\DeclareMathOperator\Ima{Im}

\begin{document}

\begin{abstract}

If $G$ is a group and $M$ a $G$-module, then we denote by $(C(G,M),d)$
the corresponding cochain complex obtained from the standard
resolution. If $n\geq 0$, then an element of $H^n(G,M)$ is written as
$[a]$ for some cocycle $a\in C^n(G,M)$.

Our first result regards the action of $G$ on $H^n(G,M)$, which is
known to be trivial, i.e. $s[a]=[a]$ for all $[a]\in H^n(G,M)$ and
$s\in G$. We prove that for all $a\in C^n(G,M)$ and $s\in G$ we have
$$sa-a=(h_sd+dh_s)(a),$$
where $h_s:C(G,M)\to C(G,M)[-1]$ is an explicit linear map.

The second result involves the commutativity of the cup product,
i.e.\\ $[a]\cup [b]=(-1)^{pq}t_*([b]\cup [a])$ for all $[a]\in H^p(G,M)$
and $[b]\in H^q(G,N)$. Here $t:N\otimes M\to M\otimes N$ denotes the
natural bijection. We will prove that for all $a\in H^p(G,M)$ and
$b\in H^q(G,N)$ we have
$$(-1)^{pq}t_*(b\cup a)-a\cup b=(hd+dh)(a\otimes b),$$
where $h:C(G,M)\otimes C(G,N)\to C(G,M\otimes N)[-1]$ is an explicit
linear map.

We will use these results as a pre-requisite in a future paper.

\end{abstract}
\maketitle

\section{Introduction}

Let $G$ be a group. For any $G$-module $M$ we denote by $(C(G,M),d)$
the corresponding cochain complex obtained from the standard
resolution. If $n\geq 0$, then an element of $H^n(G,M)$ is written as
the class $[a]$ of some cocycle $a\in C^n(G,M)$.

We consider two results regarding the cohomology, namely the
triviality of the action of $G$ on $H^n(G,M)$, i.e. $s[a]=[a]$ for
every $[a]\in H^n(G,M)$ and $s\in G$, and the commutativity of the cup
product, i.e.  $[a]\cup [b]=(-1)^{pq}t_*([b]\cup [a])$ for every
$[a]\in H^p(G,M)$ and $[b]\in H^q(G,N)$. We will show how these
formulas translate to cochains. Namely, we will write the differences
$sa-a$ and $(-1)^{pq}t_*(b\cup a)-a\cup b$ in terms of some
explicit homotopies. Our main results are Theorems 1.1, 3.4 and 3.8.

We now give a short review of notations and results we will use. 

If $f:M'\to M$ and $g:N\to N'$ are morphisms of $R$-modules, then we
denote by $f^*:Hom_R(M,N)\to Hom_R(M',N)$ and
$g_*:Hom_R(M,N)\to Hom_R(M,N')$ the induced morphisms, $f^*(a)=af$ and
$g_*(a)=ga$. Note that $f^*g_*=g_*f^*$, as both maps are given by
$a\mapsto fag$.

If $(F,\partial )$ is a chain complex of $R$-modules, with
$\partial :F\to F[-1]$ the boundary map, then for every $R$-module
$M$ we have the cochain complex $(Hom_R(F,M),d)$, where
$Hom_R(F,M)^n=Hom_R(F_n,M)$ and $d:Hom_R(F,M)\to Hom_R(F,M)[1]$ is
given by $da(x)=a(\partial x)$, i.e. $d=\partial^*$.

If $(F,\partial )$ and $(F',\partial )$ are two chain complexes, 
and $f,g:F\to F'$ are two maps of complexes, we denote by
$h:f\approx g$ the fact that $h:F\to F'[1]$ is a homotopy from $f$ to
$g$, i.e. $f-g=\partial h+h\partial$. Similarly, if $(P,d)$ and
$(P',d)$ are cochain complexes and $f,g:P\to P'$ are two maps of
complexes, we denote by $h:f\approx g$ the fact that $h:P\to P'[-1]$
is a homotopy from $f$ to $g$, i.e. $f-g=hd+dh$.

If $f,g:F\to F'$ are maps of chain complexes of $R$-modules, and
$h:f\approx g$, then we have the maps of cochain complexes
$f^*,g^*:Hom_R(F',M)\to Hom_R(F,M)$ and the map
$h^*:Hom_R(F',M)\to Hom_R(F,M)[-1]$ satisfies $h^*:f^*\approx g^*$. 

Moreover, if $f,g$ are maps of chain or cochain complexes and
$h:f\approx g$, then $hk:fk\approx gk$ and
$\ell h:\ell f\approx\ell g$ for all compatible maps of complexes
$k,\ell$.

From now on, we fix the group $G$ and we denote by
$\varepsilon :F\to\ZZ$ the corresponding standard resolution, with
$F_n=\ZZ G^{n+1}$. We denote by $\partial :F\to F[-1]$ the boundary
map, $\partial (s_0,\ldots,s_n)
=\sum_{i=0}^n(-1)^i(s_0,\ldots,\hat s_i,\ldots,s_n)$.  Then for every
$G$-module $M$ the cochain complex $C(G,M)$ is defined as
$C(G,M)=Hom_G(F,M)$.

The action of $G$ on $C(G,M)$ is given by
$$(sa)(s_0,\ldots,s_n)=sa(s^{-1}s_0s,\ldots,s^{-1}s_ns)=a(s_0s,\ldots,s_ns)$$
(For inhomogeneous cochains we have
$(sa)(s_1,\ldots,s_n)=sa(s^{-1}s_1s,\ldots,s^{-1}s_ns)$.)

For every $s\in G$ we have an augmentation-preserving chain map
$\tau_s:F\to F$ given by
$(s_0,\ldots,s_n)\mapsto (s_0s,\ldots,s_ns)$. Hence for
$a\in C(G,M)=Hom_G(F,M)$ we have $(sa)(x)=a(\tau_s(x))$ so the map  
$s:Hom_G(F,M)\to Hom_G(F,M)$, i.e. $s:C(G,M)\to C(G,M)$, can be
written as $s=\tau_s^*$.
\medskip

The cup product,
$\cup :Hom_G(F,M)\otimes Hom_G(F,N)\to Hom_G(F,M\otimes N)$,
i.e. $\cup :C(G,M)\otimes C(G,N)\to C(G,M\otimes N)$, is defined as
the composition
$$Hom_G(F,M)\otimes Hom_G(F,N)\xrightarrow\times
Hom_G(F\otimes F,M\otimes N)\xrightarrow{\Delta^*}Hom_G(F,M\otimes N),$$
where $\times$ is given by $(a\times b)(x\otimes y)=a(x)\otimes b(y)$
and $\Delta :F\to F\otimes F$ is the Alexander-Whitney map,
given by $(s_0,\ldots,s_n)\mapsto\sum_{i=0}^n(s_0,\ldots,s_i)\otimes
(s_i,\ldots,s_n)$,  which is known to be a diagonal approximation,
i.e. an augmentation-preserving chain map from $F$ to $F\otimes F$.
See, e.g., [E, \S 3.1]. 

(By abuse of notation, here $\times$ and $\cup$ denote both two binary
operations on $Hom_G(F,M)\times Hom_G(F,N)$ and two maps on
$Hom_G(F,M)\otimes Hom_G(F,N)$. If $a\in Hom_G(F,M)$ and
$b\in Hom_G(F,N)$, then $\times (a\otimes b)=a\times b$ and
$\cup (a\otimes b)=a\cup b$.)

If $a\in Hom_G(F,M)^p$, $b\in Hom_G(F,N)^q$, then $d\times (a\otimes
b)=d(a\times b)=da\times b+(-1)^pa\times db=\times (da\otimes
b+(-1)^pa\otimes db)=\times d(a\otimes b)$. Hence $d\times =\times d$,
i.e. $\times$ is a map of complexes.

Following [E, \S 3.1 (d)] we consider the augmentation-preserving
chain map $\tau :F\otimes F\to F\otimes F$, given by
$x\otimes y\mapsto (-1)^{\deg x\deg y}y\otimes x$, the map of
cochain complexes
$T:Hom_R(F,M)\otimes Hom_R(F,N)\to Hom_R(F,N)\otimes Hom_R(F,M)$,
given by $a\otimes b\mapsto (-1)^{\deg a\deg b}b\otimes a$, and the
map $t:N\otimes M\to M\otimes N$, given by
$\beta\otimes\alpha\to\alpha\otimes\beta$.  

Let $\bar\cup: Hom_G(F,M)\otimes Hom_G(F,N)\to Hom_G(F,M\otimes N)$,
$\bar\cup =t_*\cup T$. If $a\in Hom_G(F,M)^p$, $b\in Hom_G(F,M)^p$,
then $a\bar\cup b=\bar\cup (a\otimes b)=t_*\cup (T(a\otimes b))=t_*\cup
((-1)^{pq}b\otimes a)=(-1)^{pq}t_*(b\cup a)$.

We note that $t_*\times T=\tau^*\times$. If $a\in Hom_G(F,M)^p$ and
$b\in Hom_G(F,M)^p$, then
$t_*\times T(a\otimes b)=(-1)^{pq}t(b\times a)$ and
$\tau^*\times (a\otimes b)=(a\times b)\tau$. If $x\in F_k$,
$y\in F_l$, then
$(-1)^{pq}t(b\times a)(x\otimes y)=(-1)^{pq}a(y)\otimes b(x)$ and
$(a\times b)\tau (x\otimes y)=(-1)^{kl}a(y)\otimes b(x)$. But
$(-1)^{pq}a(y)\otimes b(x)=(-1)^{kl}a(y)\otimes b(x)$. (If
$(k,l)=(q,p)$, this is obivious. Otherwise both sides are zero.)
\footnote{The formula appears in the equivalent form
$\times T=t_*\tau^*\times$ in [E,\S 3.1(d)], but on cohomology instead
of cochains.}

We also have $t_*\Delta^*=\Delta^*t_*$. Hence
$\bar\cup =t_*\cup T=t_*\Delta^*\times T=\Delta^*t_*\times
T=\Delta^*\tau^*\times$.

\btm (i) If $s\in G$ and $\phi_s:F\to F[1]$ satisfies
$\phi_s:\tau_s\approx 1$, then for every $G$-module $M$, the map
$h_s=\phi_s^*:C(G,M)\to C(G,M)[-1]$, given by $h_s(a)=a\phi_s$, satisfies
$h_s:s\approx 1$, i.e.
$$sa-a=(h_sd+dh_s)(a)~\forall a\in C(G,M).$$

(ii) If $\lambda :F\to (F\otimes F)[1]$ satisfies
$\lambda :\tau\Delta\approx\Delta$, then for every $G$-modules $M$ and
$N$ the map
$h=\lambda^*\times:C(G,M)\otimes C(G,N)\to C(G,M\otimes N)[-1]$, given
by $h(a\otimes b)=(a\times b)\lambda$, satisfies
$h:\bar\cup\approx\cup$, i.e.
$$(-1)^{pq}t_*(b\cup a)-a\cup b=(hd+dh)(a\otimes b)~\forall a\in
C^p(G,M),\, b\in C^q(G,N).$$
\etm
\pf (i) The map $s:C(G,M)\to C(G,M)$ writes as $s=\tau_s^*$. Then,
since $\phi_s:\tau_s\approx 1$, we have $\phi_s^*:\tau_s^*\approx 1$,
i.e. $h_s:s\approx 1$.  

(ii) We have $\cup =\Delta^*\times$ and
$\bar\cup =\Delta^*\tau^*\times =(\tau\Delta )^*\times$. Since
$\lambda :\tau\Delta\approx\Delta$, we have
$\lambda^*:(\tau\Delta )^*\approx\Delta^*$. Since $\times$ is a map of
cochain complexes, we get
$\lambda^*\times :(\tau\Delta )^*\times\approx\Delta^*\times$,
i.e. $h:\bar\cup\approx\cup$.

If $a\in C^p(G,M)$ and $b\in C^q(G,N)$ then
$h(a\otimes b)=\lambda^*\times (a\otimes b)=\lambda^*(a\times
b)=(a\times b)\lambda$, as claimed. And, since
$\bar\cup (a\otimes b)=(-1)^{pq}t_*(b\cup a)$ and
$\cup (a\otimes b)=a\cup b$, the relation $\bar\cup -\cup =hd+dh$
applied to $a\otimes b$ writes as
$(-1)^{pq}t_*(b\cup a)-a\cup b=(hd+dh)(a\otimes b)$. \qed

Since both $\tau_s,\, 1:F\to F$ are augmentation-preserving chain
maps, we have $\tau_s\approx 1$ by [B, I.7.5]. Similarly,
$\tau,\, 1:F\otimes F\to F\otimes F$ are augmentation-preserving chain
maps, so $\tau\approx 1$. Since $\Delta$ is a chain map, this
implies that $\tau\Delta\approx\Delta$. (Alternatively, both
$\tau\Delta,\,\Delta :F\to F\otimes F$ are augmentation-preserving
chain maps, so $\tau\Delta\approx\Delta$.) To produce explicit maps
$h_s$ and $h$ we need explicit homotopies $\phi_s:\tau_s\approx 1$ and
$\lambda :\tau\Delta\approx\Delta$.

\section{Explicit contracting homotopies}

We start by giving an explicit formula for a contracting homotopy of a
product of two resolutions in terms of the contracting homotopies for
the two resolutions.

\blm If $h$ and $h'$ are contracting homotopies for the resolutions
$F\xrightarrow\varepsilon M\to 0$ and
$F'\xrightarrow{\varepsilon'}N\to 0$, then we have a contracting
homotopy $H$ of
$F\otimes F'\xrightarrow{\varepsilon\otimes\varepsilon'}M\otimes N\to 0$,
where $H_{-1}(x\otimes y)=h_{-1}(x)\otimes h_{-1}(y)$ for every
$x\otimes y\in M\otimes N$ and, if $x\otimes y\in F\otimes F'$, then
$$H(x\otimes y)=\begin{cases}h(x)\otimes y&\text{if }x\in F_{>0}\\
h(x)\otimes y+h_{-1}\varepsilon (x)\otimes h'(y)&\text{if }x\in F_0
\end{cases}$$
\elm
\pf This is essentially [AM, Lemma IV.3.2]. Note that in the
terminology of [AM, \S IV.3] by a contracting homotopy of
$F\xrightarrow\varepsilon M\to 0$ is meant only the restriction
$\tilde h:F\to F[1]$ of $h$ i.e. $\tilde h=(h_p)_{p\geq 0}$. We
consider the similar restrictions $\tilde h'$ of $h'$ and $\tilde H$
of $H$. Then, in the notations of [AM, Lemma IV.3.2], $\tilde h$,
$\tilde h'$ and $\tilde H$ are denoted by $s_G$, $s_H$ and
$s^\otimes$, both $\varepsilon$ and $\varepsilon'$ are denoted by
$\epsilon$ and both $h_{-1}$ and $h'_{-1}$ by $\phi$. The map
$H_{-1}$, although not explicitly stated, is $\phi\otimes\phi$, so we
have $H_{-1}=h_{-1}\otimes h'_{-1}$, i.e.
$H_{-1}(x\otimes y)=h_{-1}(x)\otimes h'_{-1}(y)$. Then
[AM, Lemma IV.3.2] writes as
$\tilde H=\tilde h\otimes 1+h_{-1}\varepsilon\otimes\tilde h'$. Here
$h_{-1}\varepsilon :F\to F$ is the extension of
$h_{-1}\varepsilon :F_0\to F_0$ so it is trivial on $F_{>0}$. Thus, if
$x\otimes y\in F\otimes F'$ then
$H(x\otimes y)=h(x)\otimes y+h_{-1}\varepsilon (x)\otimes h'(y)$ if
$x\in F_0$ and $=h(x)\otimes y$ if $x\in F_{>0}$. \qed


For the rest of the section, $F\xrightarrow\varepsilon\ZZ\to 0$ is the
standard resolution corresponding to a group $G$. We also denote by
$\bar F\xrightarrow\varepsilon\ZZ$ the normalized standard resolution,
i.e. $\bar F=F/D$, where $D\subseteq F$ is the subcomplex generated by
$(s_0,\ldots,s_n)$, with $s_i=s_{i+1}$ for some $i$.

We now produce contracting homotopies for
$F\xrightarrow\varepsilon\ZZ\to 0$ and
$F\otimes F\xrightarrow{\varepsilon\otimes\varepsilon}\ZZ\to 0$.

\blm For every $t\in G$ we have a contracting homotopy $\psi_t$ of
$F\xrightarrow\varepsilon\ZZ\to 0$, regarded as a complex of
$\ZZ$-modules, which on $\ZZ$ is given by $\psi_t(1)=t$ and on $F$ by
$\psi_t(x)=(t,x)$ for every $x$ in the basis $G^{n+1}$ of $F_n$, with
$n\geq 0$.

The contracting homotopy $\psi_t$ can also be defined on the
normalized resolution $\bar F\xrightarrow\varepsilon\ZZ\to 0$.
\elm
\pf Straightforward. For the case $t=1$, see, say, [B, \S I.3.], but
the proof is the same for $t$ arbitrary. \qed

\blm For every $s,t\in G$ we have a contracting homotopy $\psi_{s,t}$
of  $F\otimes F\xrightarrow{\varepsilon\otimes\varepsilon}\ZZ\to 0$,
regarded as a complex of $\ZZ$-modules, which on $\ZZ$ is given by
$\psi_{s,t}(1)=s\otimes t$ and on $F\otimes F$, for every
$x\otimes y$, with $x$ and $y$ in the bases $G^{p+1}$ and $G^{q+1}$ of 
$F_q$ and $F_q$, by
$$\psi_{s,t}(x\otimes y)=\begin{cases}(s,x)\otimes y&
\text{if }p>0\\
(s,x)\otimes y+s\otimes (t,y) &\text{if }p=0\end{cases}.$$

Moreover, $\psi_{s,t}$ can also be defined on the normalized
resolution, i.e. on\\ 
$\bar F\otimes\bar F\xrightarrow{\varepsilon\otimes\varepsilon}\ZZ\to 0$.  
\elm
\pf We apply Lemma 2.1 where both $F\xrightarrow\varepsilon M\to 0$
and $F'\xrightarrow{\varepsilon'} N\to 0$ are 
$F\xrightarrow\varepsilon\ZZ\to 0$, and the contracting homotopies $h$
and $h'$ are $\phi_s$ and $\phi_t$, as defined in Lemma 2.2. We prove
that the resulting contracting homotopy $H$ of
$F\otimes F\xrightarrow{\varepsilon\otimes\varepsilon}\ZZ\to 0$ is
$\psi_{s,t}$.

For $H_{-1}:\ZZ\cong\ZZ\otimes\ZZ\to (F\otimes F)_0$, we have
$H_{-1}(1)=H_{-1}(1\otimes 1)=\psi_s(1)\otimes\psi_t(1)=s\otimes t$.
Hence $H$ coincides with $\psi_{s,t}$ on $\ZZ$.

Suppose now that $x\in G^{p+1}\subseteq F_p$ and
$y\in G^{q+1}\subseteq F_q$. If $p=0$, then $H(x\otimes
y)=\psi_s(x)\otimes y+\psi_s\varepsilon
(x)\otimes\psi_t(y)=(s,x)\otimes y+s\otimes (t,y)=\psi_{s,t}(x\otimes
y)$. (We have $\psi_s\varepsilon (x)=\psi_s(1)=s$.) If $p>0$, then
$H(x\otimes y)=\psi_s(x)\otimes y=(s,x)\otimes y=\psi_{s,t}(x\otimes
y)$, which concludes the proof.

Since $\psi_{s,t}$ is given in terms of $\psi_s$ and $\psi_t$, which
can also be defined on normalized resolutions, the above reasoning can
be repeated with $F$ replaced by $\bar F$, so our result also holds on
the normalized resolution. \qed

\blm On normalized resolutions we have $\psi_t^2=0$ $\forall t\in G$
and $\psi_{s,t}^2=0$ $\forall s,t\in G$.
\elm
\pf Recal that in $\bar F$ we have $(s_0,\ldots,s_n)=0$ if
$s_i=s_{i+1}$ for some $i$.

We have $\psi_t^2(1)=\psi_t(t)=(t,t)=0$ so $\psi_t^2=0$ on
$\ZZ$. For $\bar F_n$ with $n\geq 0$ we note that if $x\in G^{n+1}$,
then $\psi^2(x)=(t,t,x)=0$.

We have $\psi_{s,t}^2(1)=\psi_{s,t}(s\otimes t)=(s,s)\otimes
t+s\otimes (t,t)=0$ so $\psi_{s,t}^2=0$ on $\ZZ$. For
$\bar F_p\otimes\bar F_q$ let
$x\otimes y\in G^{p+1}\otimes G^{q+1}$. If $p>1$, then
$\psi_{s,t}^2(x\otimes y)=(s,s,x)\otimes y=0$. If $p=0$, then
$\psi_{s,t}^2(x\otimes y)=\psi_{s,t}((s,x)\otimes y+s\otimes
(t,y))=(s,s,x)\otimes y+(s,s)\otimes (t,y)+s\otimes (t,t,y)=0$. \qed

\blm (i) If $t,u\in G$, then $\psi_{ut}u=u\psi_t$. 

(ii) If $s,t,u\in G$, $\psi_{us,ut}u=u\psi_{s,t}$.
\elm
\pf (i) On $\ZZ$ we have
$\psi_{ut}u(1)=\psi_{ut}(1)=ut=u\psi_t(1)$. On $F$, if $n\geq 0$ and
$x$ belongs to the basis $G^{n+1}$ of $F_n$, then
$\psi_{ut}(ux)=(ut,ux)=u\psi_t(x)$.

(ii) On $\ZZ$ we have
$\psi_{us,ut}u(1)=\psi_{us,ut}(1)=us\otimes ut=u\psi_{s,t}(1)$. On
$F\otimes F$, if $p,q\geq 0$ and $x\otimes y$ belongs the basis
$G^{p+1}\otimes G^{q+1}$ of $F_p\otimes F_q$, then
$\psi_{us,ut}u(x\otimes y)=\psi_{us,ut}(ux\otimes uy)$. If $p>0$, then
$\psi_{us,ut}(ux\otimes uy)=(us,ux)\otimes uy=u\psi_{s,t}(x\otimes
y)$. If $p=0$, then $\psi_{us,ut}(ux\otimes uy)=(us,ux)\otimes
uy+us\otimes (ut,uy)=u\psi_{s,t}(x\otimes y)$. \qed

\blm (i) If $n\geq 0$ and $\gamma :F_n\to F$ is a $G$-linear map such
that $\Ima\gamma\subseteq\partial F$, then for any $0\leq i\leq n$ the
map $\omega :F_n\to F$, given by
$\omega (s_0,\ldots,s_n)=\psi_{s_i}\gamma (s_0,\ldots,s_n)$, is
$G$-linear and satisfies $\partial\omega =\gamma$.

(ii) If $p,q\geq 0$ and $\gamma :F_p\otimes F_q\to F\otimes F$ is a
$G$-linear map such that $\Ima\gamma\subseteq\partial (F\otimes F)$,
then for any $0\leq i,j\leq p+q+1$ the map
$\omega :F_p\otimes F_q\to F\otimes F$, given by
$\omega ((s_0,\ldots,s_p)\otimes (s_{p+1},\ldots,s_{p+q+1}))
=\psi_{s_i,s_j}\gamma ((s_0,\ldots,s_p)\otimes
(s_{p+1},\ldots,s_{p+q+1}))$, is $G$-linear and satisfies
$\partial\omega =\gamma$.
\elm
\pf The $G$-linearity, both for (i) and (ii), follows from Lemma
2.5. Let $u\in G$. For (i), if $x=(s_0,\ldots,s_n)$, since the $i$
entry of $ux=(us_0,\ldots,us_n)$ is $us_i$, we have
$\omega (ux)=\psi_{us_i}\gamma (ux)=\psi_{us_i}u\gamma (x)$. Since, by
Lemma 2.5(i), $\psi_{us_i}u=u\psi_{s_i}$, we get
$\omega (ux)=u\psi_{s_i}\gamma (x)=u\omega (x)$. Hence $\omega$ is
$G$-linear. Similarly for (ii), if
$x=(s_0,\ldots,s_p)\otimes (s_{p+1},\ldots,s_{p+q+1})$, since the $i$
and $j$ entries of
$ux=(us_0,\ldots,us_p)\otimes (us_{p+1},\ldots,us_{p+q+1})$ are $us_i$
and $us_j$, we have
$\omega (ux)=\psi_{us_i,us_j}\gamma (ux)=\psi_{us_i,us_j}u\gamma (x)$.
Since, by Lemma 2.5(ii), $\psi_{us_i,us_j}u=u\psi_{s_i,s_j}$, we get
$\omega (ux)=u\psi_{s_i,s_j}\gamma (x)=u\omega (x)$.

For the formula $\gamma =\partial\omega$, we note that if
$(F',\partial')$ is a chain complex and $h:F'\to F'[1]$ is a
contracting homotopy, then for every $y\in\partial'F'$ we have
$y=\partial'h(y)$. Indeed, we have $\partial'y\in\partial'^2F'=0$ and
so $y=\partial'h(y)+h(\partial'y)=\partial'h(y)$. For (i) we apply
this result in the case when $F'$ is
$F\xrightarrow\varepsilon\ZZ\to 0$, $h=\psi_{s_i}$ and $y=\gamma (x)$,
with $x\in F$. Since $\gamma (x)\in\Ima\gamma\subseteq\partial F$, we
have $\gamma (x)=\partial\psi_{s_i}\gamma (x)=\partial\omega (x)$.
Hence $\gamma =\partial\omega$. The proof for (ii) is similar, but
this time $F'$ is
$F\otimes F\xrightarrow{\varepsilon\otimes\varepsilon}\ZZ\to 0$,
$h=\psi_{s_i,s_j}$ and $y=\gamma (x)$, with $x\in F\otimes F$. \qed

\section{Main results}

We now determine the homotopies $\phi_s:\tau_s\approx 1$ and
$\lambda :\tau\Delta\approx\Delta$, from which we obtain the
homotopies $h_s$ and $h$ of Theorem 1.1. As a first step for
determining $\lambda :\tau\Delta\approx\Delta$, we first find a
homotopy $\phi :\tau\approx 1$. 

\bff We will use the standard procedure for constructing homotopies
from, say, [B, Lemma I.7.4], applied to augmentation-preserving maps,
as in [B, Theorem I.7.5]. If $F\xrightarrow\varepsilon M\to 0$ and
$F'\xrightarrow{\varepsilon'}M\to 0$ are projective resolutions and
$f,g:F\to F'$ are augmentation-preserving maps, then the homotopy
$h=(h_n)_{n\in\ZZ}$ is constructed inductively as follows. For $n<0$
we have $h_n=0$. For $n=0$ we have
$\varepsilon'f_0=\varepsilon'g_0=\varepsilon$, so
$\varepsilon'(f_0-g_0)=0$, which implies that
$$\Ima (f_0-g_0)\subseteq\ker\varepsilon'=\Ima (\partial':F'_1\to F'_0),$$
from which we deduce the existence of some $h_0:F_0\to F'_1$ with
$f_0-g_0=\partial'h_0$.

For the inductive step $n-1\to n$, we have
$\partial'(f_n-g_n)=(f_{n-1}-g_{n-1})\partial
=(\partial'h_{n-1}+h_{n-2}\partial )\partial =\partial'h_{n-1}\partial$,
so $\partial'(f_n-g_n-h_{n-1}\partial )=0$, which implies that
$$\Ima (f_n-g_n-h_{n-1}\partial )\subseteq\ker (\partial':F'_n\to
F'_{n-1})=\Ima (\partial':F'_{n+1}\to F'_n),$$
from which we deduce the existence of some $h_n:F_n\to F'_{n+1}$ with
$f_n-g_n-h_{n-1}\partial =\partial'h_n$,
i.e. $f_n-g_n=\partial'h_n+h_{n-1}\partial$.

In our case, when the maps $f,g:F\to F'$ are $\tau_s,1:F\to F$ or
$\tau,1:F\otimes F\to F\otimes F$, we use Lemma 2.6 to produce $h_0$
with $f_0-g_0=\partial'h_0$ or, for $n>0$, $h_n$ with
$f_n-g_n-h_{n-1}\partial =\partial'h_n$.
\eff

\bff Recall the notation $\bar F=F/D$ for the normalized resolution
from \S2, where $D\subseteq F$ is the  degenerate subcomplex generated
by $(s_0,\ldots,s_n)\in G^{n+1}$ with $s_i=s_{i+1}$ for some $i$.

The maps $\partial :F\to F[-1]$ and $\psi_t:F\to F[1]$, with $t\in G$,
can be defined also on $\bar F=F/D$, so $\partial D\subseteq D$ and
$\psi_t(D)\subseteq D$. Also note that $\tau_s(D)=D$, so
$\tau_s:F\to F$ can also be defined as $\tau_s:\bar F\to\bar F$. 

We have $F=D\oplus F'$, where $F'$ is the $\ZZ$-submodule of $F$
generated by $(s_0,\ldots,s_n)\in G^{n+1}$ with $s_i\neq s_{i+1}$ for all
$i$. 

We also denote by $F''$ the $\ZZ$-submodule of $F$
generated by $(s_0,\ldots,s_n)\in G^{n+1}$ with $s_0,\ldots,s_n$
mutually distinct and by $F''_s$ the submodule of $F$ generated by
$(s_0,\ldots,s_n)$ with $s_0,\ldots,s_n,s_0s,\ldots,s_ns$ mutually
distinct. Obviously $F''_s\subseteq F''\subseteq F'$.  If
$x=(s_0,\ldots,s_n)$ is one of the generators of $F''$ or $F''_s$,
then so is $(s_0,\ldots,\hat s_i,\ldots,s_n)$ $\forall i$. Then
$\partial x=\sum_{i=1}^n(-1)^i(s_0,\ldots,\hat s_i,\ldots,s_n)\in F''$
or $F''_s$, respectively. Thus $\partial F''\subseteq F''$ and
$\partial F''_s\subseteq F''_s$.

Also note that $\tau_s$ sends bijectively the set of generators of
$F''$ to itself, so $\tau_s(F'')=F''$.
\eff

\bpr We have $\phi_s:\tau_s\approx 1$, where $\phi_s:F\to F[1]$ is
given by
$(s_0,\ldots,s_n)\mapsto\sum_{i=0}^n(-1)^i(s_0,\ldots,s_i,s_is,\ldots,s_ns)$.
\epr
\pf We use the standard procedure discribed in 3.1 to define
inductively on $F_n$, for $n\geq 0$, the map
$\tilde\phi_s=(\tilde\phi_{s,n}):F\to F[1]$ such that
$\phi_s:\tau_s\approx 1$. 

By 3.1, on $F_0$ we have $\Ima (\tau_s-1)\subseteq\partial F$.
Then, by  Lemma 2.6(i), the map $\tilde\phi_{s,0}:F_0\to F_1$, given
by $\tilde\phi_{s,0}(s_0)=\psi_{s_0}(\tau_s-1)(s_0)$, is $G$-linear
and satisfies the desired relation
$\tau_s-1=\partial\tilde\phi_{s,0}$. 

Assume now that $n\geq 1$ and we defined
$\tilde\phi_{s,0},\ldots,\tilde\phi_{s,n-1}$. By 3.1, the image of the
map $\tau_s-1-\tilde\phi_{s,n-1}\partial :F_n\to F$ is contained in
$\partial F$. Then, by Lemma 2.6(i), the map
$\tilde\phi_{s,n}:F_n\to F_{n+1}$, given by
$\tilde\phi_{s,n}(s_0,\ldots,s_n)=
\psi_{s_0}(\tau_s-1-\tilde\phi_{s,n-1}\partial )(s_0,\ldots,s_n)$, is 
$G$-linear and satisfies the desired relation
$\tau_s-1-\tilde\phi_{s,n-1}\partial =\partial\tilde\phi_{s,n}$.

Note that $\tilde\phi_s$ is defined recursively in terms of the maps
$\tau_s$, $\partial$ and $\psi_t$, for some $t\in G$. Since all these
maps can also be defined on the normalized cochains $\bar F$,
$\tilde\phi_s$ too can be defined on $\bar F$.

We now prove, by induction on $n$, that on $\bar F$ we have
$\tilde\phi_s=\phi_s$. So, from now on, we work on $\bar F$, where
$(s_0,\ldots,s_n)=0$ if $s_i=s_{i+1}$ for some $i$.

If $n=0$, then $\tilde\phi_s(s_0)=\psi_{s_0}(\tau_s-1)(s_0)
=\psi_{s_0}(s_0s-s_0)=(s_0,s_0s)-(s_0,s_0)=(s_0,s_0s)=\phi_s(s_0)$,
as claimed. 

Before proving the induction step, we note that, by definition, we
have $\tilde\phi_s(s_0,\ldots,s_n)\in\Ima\psi_{s_0}$ so
$\psi_{s_0}\tilde\phi_s(s_0,\ldots,s_n)\in\Ima\psi_{s_0}^2$. But, by
Lemma 2.4, we have $\psi_{s_0}^2=0$. Hence
$$\psi_{s_0}\tilde\phi_s(s_0,\ldots,s_n)=0.$$

We now prove the induction step $n-1\to n$. We have
$$\begin{aligned}\tilde\phi_s(s_0,\ldots,s_n) &
=\psi_{s_0}(\tau_s-1-\tilde\phi_s\partial )(s_0,\ldots,s_n)\\
& =\psi_{s_0}(s_0s,\ldots,s_ns)-\psi_{s_0}(s_0,\ldots,s_n)
-\psi_{s_0}\tilde\phi_s(\partial (s_0,\ldots,s_n)).\end{aligned}$$
The first two terms are $(s_0,s_0s,\ldots,s_ns)$ and
$(s_0,s_0,\ldots,s_n)=0$. For the third term, note that, with
exception of $(s_1,\ldots,s_n)$, all terms of
$\partial (s_0,\ldots,s_n)$ are of the type
$\pm (s_0,s'_1,\ldots,s'_{n-1})$. But, as seen above,
$\psi_{s_0}\tilde\phi_s(s_0,s'_1\ldots,s'_{n-1})=0$, so the
contribution of these terms to
$\psi_{s_0}\tilde\phi_s(\partial (s_0,\ldots,s_n))$ is
zero. Hence $\psi_{s_0}\tilde\phi_s(\partial
(s_0,\ldots,s_n))=\psi_{s_0}\tilde\phi_s(s_1,\ldots,s_n)$. By the
induction hypothesis, on $F_{n-1}$ we have $\tilde\phi_s=\phi_s$, so
$$\begin{aligned}
\psi_{s_0}\tilde\phi_s(s_1,\ldots,s_n)&
=\psi_{s_0}\phi_s(s_1,\ldots,s_n)=\psi_{s_0}
\biggl(\sum_{i=1}^n(-1)^{i-1}(s_1,\ldots,s_i,s_is,\ldots,s_ns)\biggr)\\
& =\sum_{i=1}^n(-1)^{i-1}(s_0,\ldots,s_i,s_is,\ldots,s_ns).
\end{aligned}.$$
In conclusion,
$$\begin{aligned}
\tilde\phi_s(s_0,\ldots,s_n) & =(s_0,s_0s,\ldots,s_ns)-0
-\sum_{i=1}^n(-1)^{i-1}(s_0,\ldots,s_i,s_is,\ldots,s_ns)\\
& =\sum_{i=0}^n(-1)^i(s_0,\ldots,s_i,s_is,\ldots,s_ns)
=\phi_s(s_0,\ldots,s_n). 
\end{aligned}$$

We now return to the original resolution $F$. We proved that over
$\bar F=F/D$ we have $\tilde\phi_s=\phi_s$, i.e. that $\tilde\phi_s$
reduced modulo $D$ coincides with $\phi'_s$. Hence
$\tilde\phi_s=\phi_s+\phi'_s$, where $\phi'_s$ takes values in $D$,
i.e. $\phi'_s:F\to D[1]$.

We want to prove that, together with
$\tau_s-1=\partial\tilde\phi_s+\tilde\phi_s\partial$, we also have
$\tau_s-1=\partial\phi_s+\phi_s\partial$.

Note that if $(s_0,\ldots,s_n)$ is a generator of $F''_s$, i.e. with
$s_0,\ldots,s_n,s_0s,\ldots,s_ns$ mutually distinct, then for every
$i$ we have $(s_0,\ldots,s_i,s_is,\ldots,s_ns)\in F''$ and so
$\phi_s(s_0,\ldots,s_n)\in F''$. Hence $\phi_s(F''_s)\subseteq F''$.
We refer to 3.2 for other properties of $D$, $F'$, $F''$ and $F''_s$. 

We first prove the relation $\tau_s-1=\partial\phi_s+\phi_s\partial$
on $F''_s$. Let $x\in F''_s$. Since
$(\tau_s-1)(x)=(\partial\tilde\phi_s+\tilde\phi_s\partial )(x)$ and
$\tilde\phi_s=\phi_s+\phi'_s$, we have
$$(\tau_s-1)(x)=(\partial\phi_s+\phi_s\partial )(x)
+(\partial\phi'_s+\phi'_s\partial )(x).$$
We prove that
$(\tau_s-1)(x),\, (\partial\phi_s+\phi_s\partial )(x)\in F'$ and
$(\partial\phi'_s+\phi'_s\partial )(x)\in D$. Then, since
$F=D\oplus F'$, the relation above implies that
$(\tau_s-1)(x)=(\partial\phi_s+\phi_s\partial )(x)$. 

We have $x\in F''_s\subseteq F'$, so $\tau_s(x)\in\tau_s(F')=F'$. Thus
$(\tau_s-1)(x)\in F'$. We have
$\phi_s(x)\in\phi_s(F''_s)\subseteq F''$ and
$\partial x\in\partial F''_s\subseteq F''_s$. It follows that
$\partial\phi_s(x)\in\partial (F'')\subseteq F''$ and
$\phi_s(\partial x)\in\phi_s(F''_s)\subseteq F''$. Hence
$(\partial\phi_s+\phi_s\partial )(x)\in F''\subseteq F'$. Finally,
since $\phi'_s$ takes  values in $D$ and $\partial D\subseteq D$, the
images of both $\partial\phi'_s$ and $\phi'_s\partial$ are included in
$D$. In particular, $(\partial\phi'_s+\phi'_s\partial )(x)\in D$. This
concludes the proof.

For the general case, denote by $t_1,\ldots,t_{2n+2}$ the sequence
$s_0,\ldots,s_n,s_0s,\ldots,s_ns$. Then, by the way $\tau_s$ and
$\phi_s$ are defined, we have 
$$(\tau_s-1-\partial\phi_s-\phi_s\partial )(s_0,\ldots,s_n)
=\sum_{1\leq i_l\leq 2n+2}
\alpha_{i_0,\ldots,i_n}(t_{i_0},\ldots,t_{i_n}),$$
for some $\alpha_{i_0,\ldots,i_n}\in\ZZ$. It suffices to prove
that all $\alpha_{i_0,\ldots,i_n}$ are zero. We consider an
arbitrary group $G$ that contains the elements $s,s_0,\ldots,s_n$ such
that $s_0,\ldots,s_n,s_0s,\ldots,s_ns$ are mutually distinct,
i.e. $t_i\neq t_j$ for $i\neq j$. (E.g. $s_i=i$ and $s=n+1$, with
$G=\ZZ$.) In this case, $x:=(s_0,\ldots,s_n)\in F''_s$, which implies
that $(\tau_s-1)(x)=(\partial\phi_s+\phi_s\partial )(x)$ and so
$\sum_{1\leq i_l\leq 2n+2}\alpha_{i_0,\ldots,i_n}
(t_{i_0},\ldots,t_{i_n})=0$. But $t_1,\ldots,t_{2n+2}$ are
mutually distinct, so $(t_{i_0},\ldots,t_{i_n})$ are mutually
distict elements in the basis $G^{n+1}$ of $F_n$. It follows that
all $\alpha_{i_0,\ldots,i_n}$ are zero. \qed

\btm (i) On $C(G,M)$ we have $h_s:\tau_s\approx 1$, where
$h_s:C^{n+1}(G,M)\to C^n(G,M)$ is given by $h_s(a)(s_0,\ldots,s_n)
=\sum_{i=0}^n(-1)^ia(s_0,\ldots,s_i,s_is,\ldots,s_ns)$.

(ii) In inhomogeneous cochains, $h_s:C^{n+1}(G,M)\to C^n(G,M)$ is
given by $h_s(a)(s_1,\ldots,s_n)=
\sum_{i=0}^n(-1)^ia(s_1,\ldots,s_i,s,s^{-1}s_{i+1}s,\ldots,s^{-1}s_ns)$.
\etm
\pf (i) Note that $h_s(a)=a\phi_s$, i.e. $h_s=\phi_s^*$, where
$\phi_s$ is the map from Proposition 3.3. Since
$\phi_s:\tau_s\approx 1$, by Theorem 1.1(i), we have
$h_s:s\approx 1$. 

(ii) Here we use the notation $C(G,M)$ for inhomogeneous cochains and
$\bar C(G,M)$ for the homogeneous ones. The homotopy $h_s$ defined on
$\bar C(G,M)$ induces a homotopy $h_s$ defined on $C(G,M)$ via the
isomorphism $C(G,M)\cong\bar C(G,M)$.

If $a\in C^n(G,M)$, we denote by $\bar a$ its correspondent in
$\bar C(G,M)$. We have $a(s_1,\ldots,a_n)=\bar a(t_0,\ldots,t_n)$,
where $t_i=s_1\cdots s_i$. (In particular, $t_0=1$.) And
$\bar a(t_0,\ldots,t_n)=t_0a(t_0^{-1}t_1,\ldots,t_{n-1}^{-1}t_n)$.

If $a\in C^{n+1}(G,M)$, $(s_1,\ldots,s_n)\in G^n$ and $t_i=s_1\cdots
s_i$, then
$h_s(a)(s_1,\ldots,s_n)=\overline{h_s(a)}(t_0,\ldots,t_n)
=\sum_{i=0}^n(-1)^i\bar a(t_0,\ldots,t_i,t_is,\ldots,t_ns)$. But
$t_0=1$, $t_{j-1}^{-1}t_j=s_j$ for $1\leq j\leq i$, $t_i^{-1}t_is=s$
and $(t_{j-1}s)^{-1}t_js=s^{-1}t_{j-1}^{-1}t_js=s^{-1}s_js$ for
$i+1\leq j\leq n$. Hence $\bar a(t_0,\ldots,t_i,t_is,\ldots,t_ns)
=a(s_1,\ldots,s_i,s,s^{-1}s_{i+1}s,\ldots,s^{-1}s_ns)$, which implies
the claimed formula for $h_s(a)(s_1,\ldots,s_n)$. \qed

\bff The map $\tau :F\otimes F\to F\otimes F$, given by
$x\otimes y\mapsto (-1)^{pq}y\otimes x$ if $x\in F_p$, $y\in F_q$ can
obviously be defined as $\bar F\otimes\bar F\to\bar F\otimes\bar F$.
Same happens for $\partial :F\otimes F\to F\otimes F[-1]$, defined as
$x\otimes y\mapsto\partial x\otimes y+(-1)^px\otimes\partial y$,
because $\partial :F\to F[-1]$ can also be defined as
$\partial :\bar F\to\bar F[-1]$. And, by Lemma 2.3, $\psi_{s,t}$ too
can be defined on the normalized resolution.

We refer to 3.2 for the notations $D$, $F'$ and $F''$.

Since $\bar F=F/D$, we have $\bar F\otimes\bar F\cong (F\otimes F)/C$,
where $C=D\otimes F+F\otimes D$.

Then, since $\partial :F\otimes F\to F\otimes F$ induces a map
$\partial :\bar F\otimes\bar F\to\bar F\otimes\bar F$, and
$\bar F\otimes\bar F=(F\times F)/C$, we have $\partial C\subseteq C$.

Since $F=D\oplus F'$, we have $F\otimes F=C\oplus (F'\otimes F')$.

We obviously have $\tau (F'\otimes F')=F'\otimes F'$.

By 3.2, $\partial F''\subseteq F''$, which implies that
$\partial (F''\otimes F'')\subseteq F''\otimes F''$. (We have
$\partial (x\otimes y)=\partial x\otimes y+(-1)^{\deg
x}x\otimes\partial y$.)  

We denote by $C'$ the $\ZZ$-submodule of $F\otimes F$ generated by
$(s_0,\ldots,s_p)\otimes (s_{p+1},\ldots,s_{p+q+1})$ where
$s_0,\ldots,s_{p+q+1}$ are mutually distinct. For each such generator
we have $(s_0,\ldots,s_p)$, $(s_{p+1},\ldots,s_{p+q+1})\in F''$, so
$C'\subseteq F''\otimes F''$. Also
$\partial ((s_0,\ldots,s_p)\otimes (s_{p+1},\ldots,s_{p+q+1}))$ is a
sum of elements of the form $\pm (s_0,\ldots,\hat s_i,\ldots
s_p)\otimes (s_{p+1},\ldots,s_{p+q+1})$ and $\pm
(s_0,\ldots,s_p)\otimes (s_{p+1},\ldots,\hat s_i,\ldots s_{p+q+1})$,
which all belong to $C'$. Hence $\partial C'\subseteq C'$.

If $x=(s_0,\ldots,s_n)$ is one of the generators of $F'$, then so are
$(s_0,\ldots,s_i)$ and $(s_i,\ldots,s_n)$ for every $i$. Thus
$\Delta (x)=\sum_{i=0}^n(s_0,\ldots,s_i)\otimes (s_i,\ldots,s_n)\in
F'\otimes F'$. Hence $\Delta (F')\subseteq F'\otimes F'$. 
\eff

\bpr We have $\phi:\tau\approx 1$, where
$\phi :F\otimes F\to F\otimes F[1]$ is given by
\begin{multline*}
(s_0,\ldots,s_p)\otimes y\mapsto
(-1)^{pq}\sum_{i=0}^p(-1)^{(q+1)i}(s_0,\ldots,s_i,y)\otimes
(s_i,\ldots,s_p)\\
-(-1)^p\sum_{i=0}^p(s_0,\ldots,s_i)\otimes (s_i,\ldots,s_p,y),
\end{multline*}
For every $(s_0,\ldots,s_p)\in G^{p+1}\subseteq F_p$ and
$y\in G^{q+1}\subseteq F_q$.
\epr
\pf We use the technique from  3.1 to define inductively
on $(F\otimes F)_n$, for $n\geq 0$, the map
$\tilde\phi =(\tilde\phi_n)_{n\geq 0}:F\otimes F\to F\otimes F[1]$
such that $\tilde\phi :\tau\approx 1$.

For convenience, if $s\in G$, we denote by $\bar\psi_s=\psi_{s,s}$.

If $n=0$, then, by 3.1, on $(F\otimes F)_0=F_0\otimes F_0$ we have
$\Ima (\tau -1)\subseteq\partial (F\otimes F)$. Then, by Lemma
2.6(ii), if we define $\tilde\phi_0:(F\otimes F)_0\to (F\otimes F)_1$
by
$$\tilde\phi_0(s_0\otimes y)=\bar\psi_{s_0}(\tau -1)(s_0\otimes y)$$
for every $s_0,y\in G\subseteq F_0$, then $\tilde\phi_0$ is $G$-linear
and satisfies the desired relation, $\tau -1=\partial\tilde\phi_0$, on
$(F\otimes F)_0$. 

For the inductive step $n-1\to n$, assume that $n>0$ and we
constructed $\tilde\phi_0,\ldots,\tilde\phi_{n-1}$. By 3.1, the image
of the map $\tau-1-\phi_{n-1}\partial: (F\otimes F)_n\to F\otimes F$
is contained in $\partial (F\otimes F)$. Then, by Lemma 3.2(ii), if
for every component $F_p\otimes F_q$, with $p+q=n$, of
$(F\otimes F)_n$ we define $\tilde\phi_n:F_p\otimes F_q\to (F\otimes
F)_{n+1}$ by
$$\tilde\phi_n((s_0,\ldots,s_p)\otimes y) =\bar\psi_{s_0}(\tau
-1-\tilde\phi_{n-1}\partial )((s_0,\ldots,s_p)\otimes y)$$
for every $(s_0,\ldots,s_p)\in G^{p+1}\subseteq F_p$ and
$y\in G^{q+1}\subseteq F_q$, then $\tilde\phi_n$ is $G$-linear and the
desired relation, $\tau -1-\tilde\phi_{n-1}\partial
=\partial\tilde\phi_n$, holds on
$(F\otimes F)_n=\bigoplus_{p+q=n}F_p\otimes F_q$. 

Note that the inductive definition of $\tilde\phi$ is given in terms
of the maps $\partial$, $\tau$ and $\bar\psi_s$, with $s\in G$. By
3.5, these maps can also be defined on the normalized resolution.
Therefore $\tilde\phi$ too induces a map on normalized resolutions,
$\tilde\phi :\bar F\otimes\bar F\to\bar F\otimes\bar F[1]$. Same as
for Proposition 3.3, we first determine $\tilde\phi$ on normalized
resolutions, where we have $(s_0,\ldots,s_n)=0$ whenever $s_i=s_{i+1}$
for some $i$.

We take first case $n=0$, i.e. when $p=q=0$ and so $s_0,y\in G$. We
have $\tilde\phi (s_0\otimes y)=\bar\psi_{s_0}(\tau -1)(s_0\otimes
y)=\bar\psi_{s_0}(y\otimes s_0-s_0\otimes y)=((s_0,y)\otimes
s_0+s_0\otimes (s_0,s_0))-((s_0,s_0)\otimes
y+s_0\otimes (s_0,y))=(s_0,y)\otimes s_0-s_0\otimes (s_0,y)$.

Before going further, we note that, by definition,
$\tilde\phi ((s_0,\ldots,s_p)\otimes y)\in\Ima\bar\psi_{s_0}$ whenever
$y\in G^{q+1}$ for some $q\geq 0$. By linearity, this holds for every
$y\in\bar F$. It follows that
$\bar\psi_{s_0}\tilde\phi ((s_0,\ldots,s_p)\otimes y)\in\Ima\bar\psi_{s_0}^2$.
But, by Lemma 2.4, $\bar\psi_{s_0}^2=0$. Hence
$$\bar\psi_{s_0}\tilde\phi ((s_0,\ldots,s_p)\otimes y)=0
\qquad\forall y\in\bar F.$$

We now consider the case $p=0$ and $q>0$. We have
$\partial (s_0\otimes y)=s_0\otimes\partial y$ and, by the above
property, $\bar\psi_{s_0}\tilde\phi (s_0\otimes\partial y)=0$. Thus
$$\begin{aligned}
\tilde\phi (s_0\otimes y)&=
\bar\psi_{s_0}(\tau -1-\tilde\phi\partial )(s_0\otimes y)
=\bar\psi_{s_0}(y\otimes s_0)-\bar\psi_{s_0}(s_0\otimes y)
-\bar\psi_{s_0}\tilde\phi (s_0\otimes\partial y)\\
&=(s_0,y)\otimes s_0-((s_0,s_0)\otimes y+s_0\otimes (s_0,y))-0
=(s_0,y)\otimes s_0-s_0\otimes (s_0,y).
\end{aligned}$$

We now prove, by induction on $p$, that for every $q\geq 0$,
$(s_0,\ldots,s_p)\in G^{p+1}\subseteq F_p$ and
$y\in G^{q+1}\subseteq F_q$ we have
\begin{multline*}
\tilde\phi ((s_0,\ldots,s_p)\otimes
y)=\sum_{i=0}^pa_{p,q}(i)(s_0,\ldots,s_i,y)\otimes (s_i,\ldots,s_p)\\
+\sum_{i=0}^pb_{p,q}(i)(s_0,\ldots,s_i)\otimes (s_i,\ldots,s_p,y)
\end{multline*}
for some integers $a_{p,q}(i)$ and $b_{p,q}(i)$. Later we will
identify $a_{p,q}(i)=(-1)^{pq+(q+1)i}$ and $b_{p,q}(i)=(-1)^{p+1}$,
which implies that $\tilde\phi =\phi$. 

If $p=0$, both when $q=0$ or $q>0$, we have
$\tilde\phi (s_0\otimes y)=(s_0,y)\otimes s_0-s_0\otimes (s_0,y)$ so the
claimed result holds, with $a_{0,q}(0)=1$ and $b_{0,q}(0)=-1$.

Now we prove the induction step $p-1\to p$ for $p\geq 1$. We use the
relation $\tilde\phi ((s_0,\ldots,s_p)\otimes y)=\bar\psi_{s_0}
(\tau -1-\tilde\phi\partial )((s_0,\ldots,s_p)\otimes y)$.

Note that $\partial ((s_0,\ldots,s_p)\otimes y)$ writes as
$\partial (s_0,\ldots,s_p)\otimes
y+(-1)^p(s_0,\ldots,s_p)\otimes\partial y$ and the term
$\partial (s_0,\ldots,s_p)\otimes y$ writes as
$(s_1,\ldots,s_p)\otimes y$ plus a sum of terms of the type $\pm
(s_0,s'_1,\ldots,s'_{p-1})\otimes y$. But we have both
$\bar\psi_{s_0}\tilde\phi ((s_0,\ldots,s_p)\otimes\partial y)=0$ and
$\bar\psi_{s_0}\tilde\phi ((s_0,s'_1,\ldots,s'_{p-1})\otimes y)=0$, so in
$\bar\psi_{s_0}\tilde\phi\partial ((s_0,\ldots,s_p)\otimes y)$ these terms
have no contribution. It follows that
$$\bar\psi_{s_0}\tilde\phi\partial ((s_0,\ldots,s_p)\otimes y)
=\bar\psi_{s_0}\tilde\phi ((s_1,\ldots,s_p)\otimes y)$$
We also have $\tau ((s_0,\ldots,s_p)\otimes y)
=(-1)^{pq}y\otimes (s_0,\ldots,s_p)$. It follows that
$$\tilde\phi ((s_0,\ldots,s_p)\otimes y)=\bar\psi_{s_0}((-1)^{pq}y\otimes
(s_0,\ldots,s_p) -(s_0,\ldots,s_p)\otimes y-\tilde\phi
((s_1,\ldots,s_p)\otimes y)).$$

If $q=0$ then $\bar\psi_{s_0}(y\otimes (s_0,\ldots,s_p))=(s_0,y)\otimes
(s_0,\ldots,s_p)+s_0\otimes (s_0,s_0,\ldots,s_p)=(s_0,y)\otimes
(s_0,\ldots,s_p)$. If $q>0$, then again $\bar\psi_{s_0}(y\otimes
(s_0,\ldots,s_p))=(s_0,y)\otimes (s_0,\ldots,s_p)$.

We also have $\bar\psi_{s_0}((s_0,\ldots,s_p)\otimes
y)=(s_0,s_0,\ldots,s_p)\otimes y=0$.

And for the term
$\bar\psi_{s_0}\tilde\phi ((s_1,\ldots,s_p)\otimes y)$, by the
induction hypthesis, we have
\begin{multline*}
\tilde\phi ((s_1,\ldots,s_p)\otimes y) =
\sum_{i=1}^pa_{p-1,q}(i-1)(s_1,\ldots,s_i,y)\otimes
(s_i,\ldots,s_p)\\
+\sum_{i=1}^pb_{p-1,q}(i-1)(s_1,\ldots,s_i)\otimes (s_i,\ldots,s_p,y).
\end{multline*}
Now $\tilde\phi( (s_1,\ldots,s_p)\otimes y)$ is a linear combination of
$x\otimes z$ with $x\in G^{k+1}$ and $z\in G^{l+1}$ for some
$k,l\geq 0$. For almost all these terms we have $k>0$ so
$\bar\psi_{s_0}(x\otimes z)=(s_0,x)\otimes z$. The only exception is
$b_{p-1,q}(0)s_1\otimes (s_1,\ldots,s_p,y)$. We have
$\bar\psi_{s_0}(s_1\otimes (s_1,\ldots,s_p,y))= (s_0,s_1)\otimes
(s_1,\ldots,s_p,y)+s_0\otimes (s_0,\ldots,s_p,y)$, so
$\bar\psi_{s_0}(b_{p-1,q}(0)s_1\otimes (s_1,\ldots,s_p,y))$
brings an extra term, $b_{p-1,q}(0)s_0\otimes (s_0,\ldots,s_p,y))$, to
$\bar\psi_{s_0}\tilde\phi ((s_1,\ldots,s_p)\otimes y)$.

In conclusion,
$$\begin{aligned}
&\mkern-36mu \tilde\phi( (s_0,\ldots,s_p)\otimes y)
=(-1)^{pq}(s_0,y)\otimes (s_0,\ldots,s_p)-0\\
\quad &-\sum_{i=1}^pa_{p-1,q}(i-1)(s_0,\ldots,s_i,y)\otimes
(s_i,\ldots,s_p)\\
\quad & -\sum_{i=1}^pb_{p-1,q}(i-1)(s_0,\ldots,s_i)\otimes
(s_i,\ldots,s_p,y) -b_{p-1,q}(0)s_0\otimes (s_0,\ldots,s_p,y)).
\end{aligned}$$
Hence $\tilde\phi( (s_0,\ldots,s_p)\otimes y)$ has the claimed form
and, moreover, we have $a_{p,q}(0)=(-1)^{pq}$,
$b_{p,q}(0)=-b_{p-1,q}(0)$ and, if $1\leq i\leq p$, then
$a_{p,q}(i)=-a_{p-1,q}(i-1)$ and $b_{p,q}(i)=-b_{p-1,q}(i-1)$. Also
recall that $a_{0,q}(0)=1$ and $b_{0,q}(0)=-1$.

We use repeatedly the formula $b_{p,q}(0)=-b_{p-1,q}(0)$ to get
$b_{p,q}(0)=(-)^pb_{0,q}(0)=(-1)^{p+1}$.

Next, from  $a_{p,q}(i)=-a_{p-1,q}(i-1)$ and
$b_{p,q}(i)=-b_{p-1,q}(i-1)$ we get that
$$a_{p,q}(i)=(-1)^ia_{p-i,q}(0)=(-1)^i(-1)^{(p-i)q}=(-1)^{pq+(q+1)i}$$
$$b_{p,q}(i)=(-1)^ib_{p-i,q}(0)=(-1)^i(-1)^{p-i+1}=(-1)^{p+1}.$$ 

It follows that on $\bar F\otimes\bar F$ we have $\tilde\phi =\phi$.
Since, by 3.5, $\bar F\otimes\bar F=(F\otimes F)/C$, it follows that
$\tilde\phi =\phi +\phi'$, where $\phi':F\otimes F\to F\otimes F[1]$
takes values in $C$. We want to prove that, together with
$\tau -1=\partial\tilde\phi +\tilde\phi\partial$, we also have
$\tau -1=\partial\phi +\phi\partial$.

Recall that $C'$ is generated by $(s_0,\ldots,s_p)\otimes y$, with
$y=(s_{p+1},\ldots,s_{p+q+1})$ such that $s_0,\ldots,s_{p+q+1}\in G$
are mutually distinct. In this case for all $i$ we have
$(s_0,\ldots,s_i,y)$, $(s_i,\ldots,s_p)$, $(s_0,\ldots,s_i)$,
$(s_i,\ldots,s_p,y)\in F''$. Then, by the definition of $\phi$, we
have $\phi ((s_0,\ldots,s_p)\otimes y)\in F''\otimes F''$. Hence 
$\phi (C')\subseteq F''\otimes F''$. For other properties involving
$F'\otimes F'$, $F''\otimes F''$, $C$ and $C'$ we refer to 3.5.

We first prove that $\tau -1=\partial\phi +\phi\partial$ holds on
$C'$. Let $x\in C'$. Since
$\tau -1=\partial\tilde\phi +\tilde\phi\partial$ and 
$\tilde\phi =\phi +\phi'$, we have
$$(\tau -1)(x)=(\partial\phi +\phi\partial )(x)
+(\partial\phi'+\phi'\partial)(x).$$
We prove that
$(\tau -1)(x),\, (\partial\phi +\phi\partial )(x)\in F'\otimes F'$
and $(\partial\phi'+\phi'\partial)(x)\in C$. Since
$F\otimes F=C\oplus (F'\otimes F')$, this will imply that
$(\tau -1)(x)=(\partial\phi +\phi\partial )(x)$.

We have $x\in C'\subseteq F'\otimes F'$ and so
$\tau (x)\in\tau (F'\otimes F')=F'\otimes F'$. Hence
$(\tau -1)(x)\in F'\otimes F'$. We have
$\phi (x)\in\phi (C')\subseteq F''\otimes F''$ and
$\partial x\in\partial C'\subseteq C'$. It follows that
$\partial\phi (x)\in\partial (F''\otimes F'')\subseteq F''\otimes F''$
and $\phi\partial x\in\phi (C')\subseteq F''\otimes F''$. Hence
$(\partial\phi +\phi\partial )(x)\in  F''\otimes F''\subseteq
F'\otimes F'$. And we have $\Ima\phi'\subseteq C$ and
$\partial C\subseteq C$, so the images of both $\partial\phi'$ and
$\phi'\partial$ are included in $C$. In particular, we get
$(\partial\phi'+\phi'\partial)(x)\in C$, which concludes our proof.

For the general case, we note that, by the way $\tau$ and $\phi$ are
defined, if $x=(s_0,\ldots,s_p)\otimes (s_{p+1},\ldots,s_{p+q+1})$,
then
$$(\tau -1-\partial\phi -\phi\partial)(x)=\sum_{k=0}^{p+q}
\sum_{0\leq i_l\leq p+q+1}\alpha_{k,i_0,\ldots,i_{p+q+1}}
(s_{i_0},\ldots,s_{i_k})\otimes (s_{i_{k+1}},\ldots,s_{i_{p+q+1}}),$$
for some $\alpha_{k,i_0,\ldots,i_{p+q+1}}\in\ZZ$. It suffices to prove
that all $\alpha_{k,i_0,\ldots,i_{p+q+1}}$ are zero.

We take $s_0,\ldots,s_{p+q+1}$ mutually distinct elements of some
group $G$, so that
$x=(s_0,\ldots,s_p)\otimes (s_{p+1},\ldots,s_{p+q+1})\in C'$. But in
this case we proved that\\
$(\tau -1)(x)=(\partial\phi +\phi\partial )(x)$ and so
$$\sum_{k=0}^{p+q}
\sum_{0\leq i_l\leq p+q+1}\alpha_{k,i_0,\ldots,i_{p+q+1}} 
(s_{i_0},\ldots,s_{i_k})\otimes (s_{i_{k+1}},\ldots,s_{i_{p+q+1}})=0.$$
But $s_0,\ldots,s_{p+q+1}$ are mutually distinct, so
$(s_{i_0},\ldots,s_{i_k})\otimes (s_{i_{k+1}},\ldots,s_{i_{p+q+1}})$
are mutually distinct elements of the basis
$\bigcup_{k+l=p+q}G^{k+1}\otimes G^{l+1}$ of $(F\otimes F)_{p+q}$.
It follows that all $\alpha_{k,i_0,\ldots,i_{p+q+1}}$ are zero. \qed

We now determine a homotopy $\lambda :\tau\Delta\approx\Delta$. Since
$\phi :\tau\approx 1$, we have $\phi\Delta :\tau\Delta\approx\Delta$.
However our $\lambda$ will not be $\phi\Delta$, but a simplified
version of $\phi\Delta$, with all degenerate terms removed.

\bpr Let $\lambda =(\lambda_n)_{n\geq 0}:F\to F\otimes F[1]$, with
$\lambda_n:F_n\to (F\otimes F)_{n+1}$,
$\lambda_n=\sum_{p+q=n+1}\lambda_{p,q}$, where
$\lambda_{p,q}:F_n\to F_p\otimes F_q$ are given by
$\lambda_{0,n+1}=\lambda_{n+1,0}=0$ and for $p,q>0$ by
$$\lambda_{p,q}(s_0,\ldots,s_n)=(-1)^{pq+q}
\sum_{i=0}^{p-1}(-1)^{(q+1)i}
(s_0,\ldots,s_i,s_{i+q},\ldots,s_n)\otimes (s_i,\ldots,s_{i+q}).$$

Then $\lambda :\tau\Delta\approx\Delta$.
\epr
\pf By Proposition 3.6, $\phi:\tau\approx 1$. Then, since $\Delta$ is
a map of chain complexes, we have
$\tilde\lambda :=\phi\Delta :\tau\Delta\approx\Delta$. By definition,
$\phi =\phi_1+\phi_2$, where
$$\phi_1((s_0,\ldots,s_k)\otimes y)=(-1)^{kl}\sum_{i=0}^k
(-1)^{(l+1)i}(s_0,\ldots,s_i,y)\otimes (s_i,\ldots,s_k),$$
$$\phi_2((s_0,\ldots,s_k)\otimes y)=
(-1)^{k+1}\sum_{i=0}^k(s_0,\ldots,s_i)\otimes (s_i,\ldots,s_k,y)$$
for every $(s_0,\ldots,s_k)\in G^{k+1}\subseteq F_k$ and
$y\in G^{l+1}\subseteq F_l$. (Here we changed $p,q$ from Proposition
3.6 to $k,l$.) Hence $\tilde\lambda =\lambda_1+\lambda_2$, where
$\lambda_1=\phi_1\Delta$ and $\lambda_2=\phi_2\Delta$.

For $\lambda_1$ we have
$$\begin{aligned}
\lambda_1(s_0,\ldots,s_n)&=\phi_1\biggl(\sum_{k+l=n}(s_0,\ldots,s_k)\otimes
(s_k,\ldots,s_n)\biggr)\\
&=\sum_{k+l=n}(-1)^{kl}\sum_{i=0}^k(-1)^{(l+1)i}
(s_0,\ldots,s_i,s_k,\ldots,s_n)\otimes (s_i,\ldots,s_k).
\end{aligned}$$
Note that $\lambda_1$ takes values in
$(F\otimes F)_{n+1}=\bigoplus_{p+q=n+1}(F_p\otimes F_q)$. We want, for
every $p,q\geq 0$ with $p+q=n+1$, to identify the $F_p\otimes F_q$
component of $\lambda_1$, which we denote by
$\tilde\lambda_{p,q}$. For every $k,l$ with $k+l=n$ and
every $0\leq i\leq k$ we have
$(s_i,\ldots,s_k)\in G^{k-i+1}\subseteq F_{k-i}$ and
$(s_0,\ldots,s_i,s_k,\ldots,s_n)\in G^{i+1+n-k+1}\subseteq F_{n+1+i-k}$.
Hence $(s_0,\ldots,s_i,s_k,\ldots,s_n)\otimes (s_i,\ldots,s_k)$
belongs to $F_p\otimes F_q=F_{n+1-q}\otimes F_q$ iff $q=k-i$. Then we
have $k=q+i$ and the condition $k+l=n=p+q-1$ writes as
$l=(p+q-1)-k=p-i-1$. We also have the conditions $0\leq i\leq k=q+i$
and $0\leq l=p-i-1$, which are equivalent to $0\leq i\leq p-1$.

The factors $(s_0,\ldots,s_i,s_k,\ldots,s_n)$ and $(s_i,\ldots,s_k)$
write as $(s_0,\ldots,s_i,s_{i+q},\ldots,s_n)$ and
$(s_i,\ldots,s_{i+q})$. And for the coefficient $(-1)^{kl+(l+1)i}$ we
note that $kl+(l+1)i=(q+i)(p-i-1)+(p-i)i=pq+2pi-q-qi-2i^2-i\equiv
pq+q+(q+1)i\pmod 2$.

In conclusion, for every $p,q\geq 0$ with $p+q=n+1$ we have
$$\tilde\lambda_{p,q}(s_0,\ldots,s_n)=
(-1)^{pq+q}\sum_{i=0}^{p-1}(-1)^{(q+1)i}
(s_0,\ldots,s_i,s_{i+q},\ldots,s_n)\otimes (s_i,\ldots,s_{i+q}).$$
Hence $\tilde\lambda_{p,q}=\lambda_{p,q}$, for all
$(p,q)\neq (n+1,0)$. (This includes the case $(p,q)=(0,n+1)$, when the
sum above is zero, so $\tilde\lambda_{0,n+1}=0=\lambda_{0,n+1}$.)
Since $\lambda_{n+1,0}=0$, we have $\lambda_1=\lambda +\lambda_3$,
where $\lambda_3=\tilde\lambda_{n+1,0}$, i.e.
$$\lambda_3(s_0,\ldots,s_n)=\sum_{i=0}^n(-1)^i
(s_0,\ldots,s_i,s_i,\ldots,s_n)\otimes s_i.$$
Since $(s_0,\ldots,s_i,s_i,\ldots,s_n)\in D$, we have
$\lambda_3(s_0,\ldots,s_n)\in D\otimes F$, so
$\Ima\lambda_3\subseteq D\otimes F$.

For $\lambda_2$ we have
$$\begin{aligned}
\lambda_2(s_0,\ldots,s_n)&
=\phi_2\biggl(\sum_{k+l=n}(s_0,\ldots,s_k)\otimes
(s_k,\ldots,s_n)\biggr)\\
&=\sum_{k+l=n}(-1)^{l+1}\sum_{i=0}^k
(s_0,\ldots,s_i)\otimes (s_i,\ldots,s_k,s_k,\ldots,s_n).
\end{aligned}$$
Since $(s_i,\ldots,s_k,s_k,\ldots,s_n)\in D$, we have
$\lambda_2(s_0,\ldots,s_n)\in F\otimes D$, so
$\Ima\lambda_2\subseteq F\otimes D$.

Since $\tilde\lambda =\lambda_1+\lambda_2$ and
$\lambda_1=\lambda +\lambda_3$, we have
$\tilde\lambda =\lambda+\lambda'$, where
$\lambda'=\lambda_2+\lambda_3$. Since
$\Ima\lambda_2\subseteq F\otimes D$ and
$\Ima\lambda_3\subseteq D\otimes F$, we have
$\Ima\lambda'\subseteq F\otimes D+D\otimes F=C$.

If $x=(s_0,\ldots,s_n)$ is one of the generators of $F''$, i.e. with
$s_j\neq s_h$ $\forall j\neq h$, then every term
$\pm (s_0,\ldots,s_i,s_{i+q},\ldots,s_n)\otimes (s_i,\ldots,s_{i+q})$
that appears in $\lambda (x)$ belongs to $F''\otimes F''$ so
$\lambda (x)\in F''\otimes F''$. Hence
$\lambda (F'')\subseteq F''\otimes F''$. For other properties
involving $F'$, $F''$, $F'\otimes F'$, $F''\otimes F''$ and $C$ we
refer to 3.2 and 3.5.

We now employ a technique which is similar to that from Propositions
3.3 and 3.6 to prove that
$\tau\Delta -\Delta =\partial\lambda +\lambda\partial$. First we prove
this relation on $F''$.

Let $x\in F''$. Since 
$\tau\Delta -\Delta =\partial\tilde\lambda +\tilde\lambda\partial$ and
$\tilde\lambda =\lambda +\lambda'$, we have 
$$(\tau\Delta -\Delta )(x)=(\partial\lambda +\lambda\partial )(x)
+(\partial\lambda'+\lambda'\partial )(x).$$
We prove that $(\tau\Delta -\Delta )(x),\, (\partial\lambda
+\lambda\partial )(x)\in F'\otimes F'$ and
$(\partial\lambda'+\lambda'\partial )(x)\in C$. Since
$F\otimes F=C\oplus (F'\otimes F')$, this will imply that
$(\tau\Delta -\Delta )(x)=(\partial\lambda +\lambda\partial )(x)$.

We have $x\in F''\subseteq F'$ so $\Delta (x)\in F'\otimes F'$ and
also $\tau\Delta (x)\in\tau (F'\otimes F')=F'\otimes F'$.
Thus $(\tau\Delta -\Delta )(x)\in F'\otimes F'$. We have
$\lambda (x)\in\lambda (F'')\subseteq F''\otimes F''$ and
$\partial x\in\partial F''\subseteq F''$, so $\partial\lambda
(x)\in\partial (F''\otimes F'')\subseteq F''\otimes F''$ and
$\lambda (\partial x)\in\lambda (F'')\subseteq F''\otimes F''$. Hence
$(\partial\lambda +\lambda\partial )(x)\in F''\otimes F'' \subseteq
F'\otimes F'$. And, since $\Ima\lambda'\subseteq C$ and
$\partial C\subseteq C$, the images of both $\partial\lambda'$ and
$\lambda'\partial$ are included in $C$. It follows that
$(\partial\lambda'+\lambda'\partial )(x)\in C$, which concludes the
proof.

For the general case, we note that, by the way $\tau$, $\Delta$ and
$\lambda$ are defined, the map
$\tau\Delta -\Delta -(\partial\lambda +\lambda\partial ):F_n\to
(F\otimes F)_n$ is given by
$$(s_0,\ldots,s_n)\mapsto\sum_{k=0}^n\sum_{0\leq i_l\leq n}
\alpha_{k;i_0,\ldots,i_{n+1}}(s_{i_0},\ldots,s_{i_k})\otimes
(s_{i_{k+1}},\ldots,s_{i_{n+1}})$$
for some $\alpha_{k;i_0,\ldots,i_{n+1}}\in\ZZ$. We will prove that all
$\alpha_{k;i_0,\ldots,i_{n+1}}$ are zero, which concludes our
proof. We take $s_0,\ldots,s_n$ mutually distinct elements of some
arbitrary group $G$. Then $x:=(s_0,\ldots,s_n)\in F''$, which implies
that
$(\tau\Delta -\Delta )(x)=(\partial\lambda +\lambda\partial )(x)$.
Hence in this case we have
$\sum_{k=0}^n\sum_{0\leq i_l\leq n}
\alpha_{k;i_0,\ldots,i_{n+1}}(s_{i_0},\ldots,s_{i_k})\otimes
(s_{i_{k+1}},\ldots,s_{i_{n+1}})=0$. But $s_0,\ldots,s_n$ are mutually
distinct, so
$(s_{i_0},\ldots,s_{i_k})\otimes (s_{i_{k+1}},\ldots,s_{i_{n+1}})$ are
mutually distinct elements in the basis
$\bigcup_{k+l=n}G^{k+1}\otimes G^{l+1}$ of $(F\otimes F)_n$. It
follows that all $\alpha_{k;i_0,\ldots,i_{n+1}}$ are zero. \qed

\btm (i) Let $M,N$ be $G$-modules and let
$h:C(G,M)\otimes C(G,N)\to C(G,M\otimes N)[-1]$, where if
$a\in C^p(G,M)$ and $b\in C^q(G,N)$, then
$h(a\otimes b)\in C^n(G,M\otimes N)$, with $n=p+q-1$, is defined by
$h(a\otimes b)=0$ if $p$ or $q=0$ and if $p,q>0$, by
$$(s_0,\ldots,s_n)\mapsto
(-1)^{pq+q}\sum_{i=0}^{p-1}(-1)^{q(i+1)}
a(s_0,\ldots,s_i,s_{i+q},\ldots,s_n)\otimes
b(s_i,\ldots,s_{i+q}).$$
Then $(-1)^{pq}t(b\cup a)-a\cup b=(hd+dh)(a\otimes b)$
$\forall a\in C(G,M)$, $b\in C(G,N)$.

(ii) In terms of inhomogeneous cochains,
$h(a\otimes b)$ for $p,q>0$ is given by
\begin{multline*}
(s_1,\ldots,s_n)\mapsto (-1)^{pq+q}\sum_{i=0}^{p-1}(-1)^{q(i+1)}
a(s_1,\ldots,s_i,s_{i+1}\cdots s_{i+q},s_{i+q+1},\ldots,s_{p+q-1})\\
\otimes s_1\cdots s_ib(s_{i+1},\ldots,s_{i+q}).
\end{multline*}
\etm
\pf (i) The map $\lambda$ from Proposition 3.7 satisfies
$\lambda :\tau\Delta\approx\Delta$. We prove that
$h(a\otimes b)=(a\times b)\lambda$. Then, by Theorem 1.1(ii), $h$
satisfies the required property.

Let $a\in C^p(G,M)$ and $b\in C^q(G,N)$ and let $n=p+q-1$,
i.e. $p+q=n+1$. Then $a\times b\in Hom(F\otimes F,M\otimes N)$, given
by $(a\times b)(x\otimes y)=a(x)\otimes b(y)$, is zero everywhere
outside $F_p\otimes F_q\subseteq (F\otimes F)_{n+1}$. Since the
$F_p\otimes F_q$ component of $\lambda$ is $\lambda_{p,q}$ we have
$(a\times b)\lambda =(a\times b)\lambda_{p,q}$. If $p$ or $q=0$
then $\lambda_{p,q}=0$ so $(a\times b)\lambda =0=h(a\otimes b)$. If
$p,q>0$, then each term
$\pm (s_0,\ldots,s_i,s_{i+q},\ldots,s_n))\otimes (s_i,\ldots,s_{i+q})$
from the formula for $\lambda_{p,q}(s_0,\ldots,s_n)$ is mapped by
$a\times b$ to $\pm a(s_0,\ldots,s_i,s_{i+q},\ldots,s_n))\otimes
b(s_i,\ldots,s_{i+q})$, which proves the claimed formula,
$(a\times b)\lambda (s_0,\ldots,s_n)=h(a\otimes b)(s_0,\ldots,s_n)$.

(ii) We use the notations from the proof of Theorem 3.4(ii). We have\\
$h(a\otimes b)(s_1,\ldots,s_n)
=\overline{h(a\otimes b)}(t_0,\ldots,t_n)
=h(\bar a\otimes\bar b)(t_0,\ldots,t_n)$. For the terms
$\pm \bar a(t_0,\ldots,t_i,t_{i+q},\ldots,t_n)\otimes
\bar b(t_i,\ldots,t_{i+q})$ which appear in
$h(\bar a\otimes\bar b)(t_0,\ldots,t_n)$, we note that
$t_{j-1}^{-1}t_j=s_j$, $t_i^{-1}t_{i+q}=s_{i+1}\cdots s_{i+q}$,
$t_0=1$ and $t_i=s_0\cdots s_i$. Then
$$\bar a(t_0,\ldots,t_i,t_{i+q},\ldots,t_n)
=a(s_1,\ldots,s_i,s_{i+1}\cdots s_{i+q},s_{i+q+1},\ldots,s_n)$$
$$\bar b(t_i,\ldots,t_{i+q})=s_0\cdots s_ib(s_{i+1},\ldots,b_{i+q}).$$
Hence the conclusion. \qed

\bff {\bf Alternative formulas.} Note that
$\tau_s\tau_{s^{-1}}=1$. Since $\phi_{s^{-1}}:\tau_{s^{-1}}\approx 1$
and $\tau_s$ is a map of complexes, we have
$\tau_s\phi_{s^{-1}}:\tau_s\tau_{s^{-1}}\approx\tau_s$,
i.e. $\tau_s\phi_{s^{-1}}:1\approx\tau_s$. So if
$\phi'_s=-\tau_s\phi_{s^{-1}}$, then $\phi'_s:\tau_s\approx 1$.

If $x=(s_0,\ldots,s_n)$, then $\phi_{s^{-1}}(x)
=\sum_{i=0}^n(-1)^i(s_0,\ldots,s_i,s_is^{-1},\ldots,s_ns^{-1})$,
so $\phi'_s(x)=-\tau_s\phi_{s^{-1}}(x)
=-\sum_{i=0}^n(-1)^i(s_0s,\ldots,s_is,s_i,\ldots,s_n)$. 

If we replace $\phi_s$ by $\phi'_s$, then the map $h_s$ from Theorem
3.4, satisfying $sa-a=(h_sd+dh_s)(a)$, will be replaced by $h'_s$,
given by $h'_s(a)=a\phi'_s$, i.e.
$$h'_s(a)(s_0,\ldots,s_n)
=-\sum_{i=0}^n(-1)^ia(s_0s,\ldots,s_is,s_i,\ldots,s_n).$$
In terms of inhomogeneous cochains
$$h'_s(a)(s_1,\ldots,s_n)=-\sum_{i=0}^n
(-1)^isa(s^{-1}s_1s,\ldots,s^{-1}s_is,s^{-1},s_{i+1},\ldots,s_n).$$
\smallskip

We have $\tau^2=1$ and $\tau$ is a map of complexes, so
$\lambda :\tau\Delta\approx\Delta$ implies that
$\tau\lambda :\tau^2\Delta\approx\tau\Delta$, i.e.
$\tau\lambda :\Delta\approx\tau\Delta$. It follows that
$\lambda':\tau\Delta\approx\Delta$, where $\lambda'=-\tau\lambda$.

If $p+q=n+1$, we denote by $\lambda'_{p,q}$ the $F_p\otimes F_q$
component of $\lambda'$. Since $\tau (F_q\otimes F_p)=F_p\otimes F_q$
and the $F_q\otimes F_p$ component of $\lambda$ is $\lambda_{q,p}$, we
have $\lambda'_{p,q}=-\tau\lambda_{q,p}$. If $p$ or $q=0$ then
$\lambda_{q,p}=0$, so $\lambda'_{p,q}=0$. If $p,q>0$ then
$$\lambda_{q,p}(s_0,\ldots,s_n)=(-1)^{qp+p}
\sum_{i=0}^{q-1}(-1)^{(p+1)i}
(s_0,\ldots,s_i,s_{i+p},\ldots,s_n)\otimes (s_i,\ldots,s_{i+p})$$
so
$$\lambda'_{p,q}(s_0,\ldots,s_n)=-(-1)^p\sum_{i=0}^{q-1}(-1)^{(p+1)i}
(s_i,\ldots,s_{i+p})\otimes (s_0,\ldots,s_i,s_{i+p},\ldots,s_n).$$

If we replace $\lambda$ by $\lambda'$, then the map $h$ from Theorem
3.8, with $(-1)^{pq}t(b\cup a)-a\cup b=(hd+dh)(a\otimes b)$, is
replced by $h'$, given by $h'(a\otimes b)=(a\times b)\lambda'$.

If $p$ or $q=0$, then $h'(a\otimes b)=0$. If $p,q>0$, then
$h'(a\otimes b)$ is given by
$$(s_0,\ldots,s_n)\mapsto -(-1)^p\sum_{i=0}^{q-1}(-1)^{(p+1)i}
a(s_i,\ldots,s_{i+p})\otimes b(s_0,\ldots,s_i,s_{i+p},\ldots,s_n).$$
In terms of inhomogeneous cochains, if $p,q>0$ then $h'(a\otimes b)$ is
given by
\begin{multline*}
(s_1,\ldots,s_n)\mapsto -(-1)^p\sum_{i=0}^{q-1}(-1)^{(p+1)i}
s_1\cdots s_i a(s_{i+1},\ldots,s_{i+p})\\
\otimes b(s_1,\ldots,s_i,s_{i+1}\cdots s_{i+p},s_{i+p+1},\ldots,s_n). 
\end{multline*}
\eff

{\bf References}
\medskip

[AM] Alejandro Adem, R. James Milgram, Cohomology of finite groups,
Grundlehren der mathematischen Wissenschaften, volume 309.

[B] Kenneth S. Brown, Cohomology of groups, Springer 1982.

[E] Leonard Evens, The cohomology of groups, Oxford Mathematical
Monographs. Oxford University Press, New York, 1991.

\end{document}